\documentclass[10pt,leqno]{article}
\usepackage{amssymb,amsbsy,amsmath,amsfonts,amssymb,amscd, mathrsfs}
\usepackage{graphicx} 

\usepackage[all]{xy} 
\usepackage{color}

\usepackage[english]{babel}
\usepackage[T1]{fontenc}
\usepackage{indentfirst}

\makeatletter
\@addtoreset{equation}{section}
\makeatother

\newtheorem{statement}{}[section]
\newtheorem{theorem}[statement]{Theorem}
\newtheorem{lemma}[statement]{Lemma}

\newtheorem{corollary}[statement]{Corollary}
 
\newtheorem{sublemma}[statement]{Sublemma}

\newcommand\C{\mathbb C}

\newcommand\T{\mathbb T}
\newcommand\D{\mathbb D}

\newcommand\e{{\rm e}}

\newcommand\eps{\varepsilon}
\newcommand\ind{{\rm 1\kern-.30em I}}
\newcommand\qed{\hfill $\square$}
\renewcommand \Re{{\mathfrak R}{\rm e}\,}
\renewcommand \Im{{\mathfrak I}{\rm m}\,}
\let\amphi=\phi
\let\phi=\varphi

\title{\bf Approximation numbers of composition operators on the Dirichlet space}
\author{\it Pascal Lef\`evre, Daniel Li, Herv\'e Queff\'elec, \\ 
\it Luis Rodr{\'\i}guez-Piazza\footnote{Supported by a Spanish research project MTM 2009-08934.}}

\date{\footnotesize \today}

\begin{document}

\maketitle 

\noindent{\bf Abstract.} \emph{We study the decay of approximation numbers of compact composition operators on the Dirichlet space. We give upper and 
lower bounds for these numbers. In particular, we improve on a result of O.~El-Fallah, K.~Kellay, M.~Shabankhah and A.~Youssfi, on the set of contact points 
with the unit circle of a compact symbolic composition operator acting on the Dirichlet space $\mathcal{D}$. We extend their results in two directions: first, the 
contact only takes place at the point $1$. Moreover, the approximation numbers of the operator can be arbitrarily sub-exponentially small.}
\medskip

\par\medskip


\section{Introduction}
  
\subsection{Organization of the paper.} 

The paper deals with composition operators. This area is widely studied nowadays, on various spaces of analytic functions (Hardy, Bergman, Dirichlet...spaces): one may read for instance the monographs \cite{Shap-livre} or \cite{CMC} to get an overview on the subject until the nineties, and \cite{ELFKELSHAYOU} or \cite{Dirichlet} for some recent results in the framework of the Dirichlet space. It seems natural to try to apply again some of the techniques used in the framework of Hardy or Bergman spaces. Nevertheless it is far from being that simple. Actually it often turns out that the Dirichlet space is one of the most difficult ``classical'' spaces to handle. For instance, a first difficulty appears at the very beginning of the theory: the composition operators are not necessarily bounded when we only require the symbol to belong to the Dirichlet space (whereas all the composition operators are bounded on the Hardy and the Bergman spaces).

The study of approximation numbers of composition operators acting on classical spaces of analytic functions (like the Hardy or the Bergman spaces) was initiated in \cite{JAT} and \cite{estimates} by the three last named authors. In the present paper, we get interested in the same kind of questions but for composition operators acting on the Dirichlet space. Some results already appear in \cite{Dirichlet} (among other things), but we focus exclusively on this topic in the sequel.

The notations and definitions are precised in the next subsection.

In section 2, we show that some similar phenomena (as in Hardy and Bergman spaces) hold in the framework of Dirichlet spaces. More precisely,  the approximation numbers of composition operators on the Dirichlet space cannot decay more rapidly than exponentially, and this speed of convergence can only be attained for symbols verifying $\|\varphi\|_\infty<1$ (see Theorem \ref{geom decay} and Theorem \ref{seul cas geom}). On the other hand, we investigate the extremal case and it turns out that $C_\phi$ may have almost geometric decay (in particular belong to all Schatten classes) and may touch the boundary of the $\D$ (see Theorem \ref{katkat}).

In section 3, we focus on composition operators whose symbol is a cusp map. It plays the same role in the theory as the lens maps in the theory of Hardy spaces. The rate of decay of its approximation numbers is given in Theorem \ref{cusp}.

At last, in section 4., we precise Theorem \ref{katkat} and prove in Theorem \ref{improvement} that the symbol (which will be the composition of a cusp map and a peak function) may belong to both the disk algebra and the Dirichlet space, and moreover meet the boundary precisely at $\{1\}$ with a level set which is any compact set with zero logarithmic capacity. 

\subsection{Notation and background.} 

We denote by $\D$ the unit open disk of the complex plane and by $A$ the normalized area measure $dx \, dy / \pi$ of $\D$. The unit circle is denoted by 
$\T = \partial \D$. \par 
A Schur function is an analytic self-map of $\D$ and the associated composition operator is defined, formally, by $C_\phi (f) = f \circ \phi$. 
The function $\phi$ is called the symbol of $C_\phi$.
\par

The Dirichlet space $\mathcal{D}$ is the space of analytic functions $f \colon \D \to \C$ such that: 
\begin{equation} 
\| f \|_{\cal D}^2 := | f (0)|^2 + \int_\D |f ' (z) |^2 \, dA (z) < + \infty \, .
\end{equation} 
If $f (z) = \sum_{n = 0}^\infty c_n z^n$, one has:
\begin{equation} 
\| f \|_{\cal D}^2 = |c_0|^2 + \sum_{n = 1}^\infty n \, |c_n|^2 \, .
\end{equation} 
Then $\| \ \|_{\cal D}$ is a norm on ${\cal D}$, making ${\cal D}$ a Hilbert space. We consider its subspace ${\cal D}_\ast$, consisting of functions 
$f \in {\cal D}$ such that $f (0) = 0$. In this paper, we call ${\cal D}_\ast$ the \emph{Dirichlet space}.  For further information on the Dirichlet space, 
the reader may see \cite{survey} or \cite{Ross}. \par 

Recall that, whereas every Schur function $\phi$ generates a bounded composition operator $C_\phi$ on the Hardy or Bergman spaces, it is no 
longer the case for the Dirichlet space (see \cite{McCluer-Shapiro}, Proposition 3.12, for instance). \par

The Bergman space ${\mathfrak B}$ is the space of analytic functions $f \colon \D \to \C$ such that: 
\begin{displaymath} 
\| f \|_{\mathfrak B}^2 := \int_\D |f (z) |^2 \, dA (z) < + \infty \, .
\end{displaymath} 
If $f (z) = \sum_{n = 0}^\infty c_n z^n$, one has $\| f \|_{\mathfrak B}^2 = \sum_{n = 0}^\infty \frac{|c_n|^2}{n + 1}$.

\medskip

We denote by $S_n$ the truncation operator: if $f (z) = \sum_{k = 1}^\infty c_k \, z^k$ is in ${\mathcal D}_\ast$, then
\begin{equation} \label{partial sum} 
(S_n f) (z) = \sum_{k = 1}^n  c_k \, z^k \, .
\end{equation} 

The Carleson window centered at $\xi \in \T$ and of size $h \in (0, 1)$ is defined as: 
\begin{equation} 
S (\xi, h) = \{ z \in \D \, ; \ |z - \xi | < h \} \, .
\end{equation} 

The notation $A \lesssim B$ (equivalently $B \gtrsim A$) means that $A \leq C\, B$ for some constant $C > 0$, and $A \approx B$ means that we have both 
$A \lesssim B$ and $B \lesssim A$.


\section{Approximation numbers in the general case} 

\subsection{Geometric decay of the approximation numbers} 

We saw in \cite{JAT}, that the approximation numbers of composition operators on the Hardy space $H^2$ as well as on the (weighted) Bergman spaces 
${\mathfrak B}_\alpha$, $\alpha > - 1$, cannot decay more rapidly than exponentially, and that this speed of convergence can only be attained for symbols 
mapping the unit disk $\D$ into a smaller disk $r \D$, with $0 < r < 1$. In this section, we see that the same phenomenon holds for the Dirichlet space. The 
proofs will be adapted from those of \cite{JAT}. 
\par\medskip

Our first result is on the geometric decay. 

\begin{theorem} \label{geom decay} 
Let $\phi \colon \D \to \D$ be an analytic self-map inducing a bounded composition operator on ${\mathcal D}_\ast$. Then, there exist positive constants 
$c', c > 0$ and $ 0 < r < 1$ such that the approximation numbers of the composition operator $C_\phi \colon {\mathcal D}_\ast \to {\mathcal D}_\ast$ satisfy:
\begin{equation} \label{july} 
\qquad \qquad c' \sqrt{n} \, \| \phi \|_\infty^n \geq a_n (C_\phi) \geq c\, r^n \, , \qquad n = 1 , 2, \ldots
\end{equation} 
\end{theorem} 

\noindent{\bf Proof.} We introduce two notations. First,  we set, for any operator $T$ on some Hilbert space $H$:
\begin{equation} \label{beta}
\beta (T) = \liminf_{n\to\infty} \big[ a_{n} (T) \big]^{1/ n} .
\end{equation}
Next, let 

\begin{displaymath}
\phi^{\#}(z) = \frac{ \phi '(z) \, (1 - | z |^2)}{1 - | \phi (z) |^2} 
\end{displaymath}
be the pseudo-hyperbolic derivative of $\phi$; we set: 
\begin{equation} \label{crochet} 
[\phi] = \sup_{z \in \D} |\phi^{\#} (z) | = \| \phi^{\#} \|_\infty . 
\end{equation}
Note that $[\phi] \leq 1$, by the Schwarz-Pick inequality. \par
\smallskip

The upper bound is easy. We may assume that $\| \phi \|_\infty < 1$. Note that, since $C_\phi$ is bounded on ${\mathcal D}_\ast$,  we have 
 $\phi^k = C_\phi (z^k) \in {\mathcal D}_\ast$. Then: 
\begin{align*}
a_{n} (C_\phi)^2 
& \leq \| C_\phi - C_\phi S_{n - 1} \|^2 
\leq \| C_\phi - C_\phi S_{n - 1} \|_{HS}^2  
= \sum_{k = n}^\infty \frac{\| \phi^k \|_{\mathcal D}^2}{k} \\ 
& = \sum_{k = n}^\infty \int_{\D} k \, | \phi^{k - 1} (z) |^2  \, | \phi ' (z) |^2 \,  dA (z) \\ 
& = \int_{\D} \sum_{k = n}^\infty k \, | \phi^{k - 1} (z) |^2  \, | \phi ' (z) |^2 \,  dA (z) \\ 
& \leq \int_{\D} \sum_{k = n}^\infty  k \, \| \phi \|_\infty^{2 k - 2}  \, | \phi ' (z) |^2 \, dA (z) \\ 
& \leq K (\phi) \,  n \, \| \phi \|_\infty^{2 n}  \, \| \phi \|_{\mathcal D}^2 \, ,
\end{align*}
(we used that $\sum_{k = n}^\infty k \, \rho^{k - 1} = \big(\sum_{k = n}^\infty \rho^k \big)' = \frac{\rho^{n - 1} [n - (n - 1) \rho]}{(1 - \rho)^2} 
\leq\frac{1}{(1 - \rho)^2} \, n\, \rho^{n - 1}$, with $\rho = \| \phi \|_\infty^2$), implying 
\begin{equation} 
a_{n} (C_\phi) \lesssim \sqrt{n} \, \| \phi \|_{\mathcal D} \, \| \phi \|_\infty^n \approx \sqrt{n} \, \| \phi \|_\infty^n 
\end{equation} 
and $\beta (C_\phi) \leq \| \phi \|_\infty$. 
\par\smallskip

For the lower bound in \eqref{july}, we shall prove that:
\begin{equation} \label{always} 
[\phi]^2 \leq \beta (C_\phi) \leq \| \phi \|_\infty \, , 
\end{equation}
which will give the result, since for each $\kappa < [\phi]^2$, there will be some constant $c_\kappa > 0$ such that 
$a_n (C_\phi) \geq c_\kappa\, \kappa^n$, $n \geq 1$.  \par
\smallskip

The inequality \eqref{always} is obtained as in the Hardy and Bergman cases in \cite{JAT}. We may assume that $C_\phi$ is compact on ${\mathcal D}_\ast$ 
(since otherwise $\beta (C_\phi) = 1$ and the result is trivial). Now, set $\amphi_{u} (z) = \frac{u - z}{1 - \overline{u} z}$, $u\in \D$.  Then, if $\phi$ is a 
symbol with $C_\phi$ compact on ${\mathcal D}_\ast$ and $a \in \D$, let $\psi = \amphi_{\phi (a)} \circ \phi \circ \amphi_a$. Note that the compactness of 
$C_\phi$ on ${\mathcal D}_\ast$ implies its compactness on ${\mathcal D}$. Hence we can write 
$C_\psi = C_{\amphi_a} \circ C_\phi \circ C_{\amphi_{\phi (a)}}$. Now, the relations 
$\psi (0) = 0$,  $\psi ' (0) = \frac{\phi ' (a) (1 - | a |^2)}{1 - |\phi (a) |^2} = \phi^{\#} (a)$ and the diagrams: 
\begin{displaymath}
\xymatrix{
{\cal D}_{\ast}  \ar[r]^{C_{\amphi_a}} &  {\cal D}  \ar[r]^ {C_\phi} & {\cal D} \ar[r]^{C_{\amphi_{\phi (a)}}} & {\cal D}_{\ast} \, ,
}
\end{displaymath}
with $\xymatrix{0 \ar[r]^{\amphi_a} &  a  \ar[r]^ {\phi} & \phi (a) \ar[r]^{\amphi_{\phi (a)}} & 0}$, show that $\psi \in {\mathcal D}_\ast$ 
and that $C_\psi$ is also compact on $\mathcal{D}_{*}$. Now we notice that, for any compact composition operator $C_\tau$ on $\mathcal{D}_{*}$, the 
solution $\sigma$ of the K\"onig  equation 
\begin{displaymath} 
\qquad \qquad \sigma \circ \tau = \tau ' (0) \, \sigma \, , \qquad \sigma (0) = 0, \  \sigma ' (0) = 1 
\end{displaymath} 
has to belong to $\mathcal{D}_{*}$ as this would be the case for any Hilbert space of analytic functions on $\D$. Hence, if $\psi ' (0) = \phi^{\#} (a) \neq 0$, 
the sequence of eigenvalues of $C_\psi$ is $\big( [\psi ' (0)]^n\big)_{n \geq 0}$. It follows from \cite{JAT}, Lemma~3.2 (which is an easy consequence 
of Weyl's inequality) that  $\beta (C_\phi) = \beta (C_\psi) \geq | \phi^{\#} (a) |$. Since this remains trivially true when $\phi^{\#} (a) = 0$, 
Theorem~\ref{geom decay} is proved. \qed 
\bigskip

Now, we shall see that the geometric decay can take place only for symbols $\phi$ such that $\|\phi \|_\infty < 1$. 

\begin{theorem} \label{seul cas geom} 
Let $C_\phi$ a bounded composition operator on ${\mathcal D}_\ast$. Then for each $r \in (0, 1)$, there exists $s = s (r) \in (0, 1)$, with 
$s (r) \to 1$ as $r \to 1$, and such that:
\begin{equation} 
\| \phi \|_\infty > r \quad \Longrightarrow \quad a_n (C_\phi) \gtrsim \frac{s^n}{\sqrt n} \, \cdot
\end{equation} 
\end{theorem} 

We shall see in the proof that we can take $s = \e^{- \eps \pi}$, with $\eps = 2 \pi /\log \frac{1 + r}{1 - r}$ (see \eqref{notation}, where $s$ is changed into 
$s^2$). \par 

Note that, in particular, with the notation \eqref{beta}, one has: 
\begin{displaymath} 
\| \phi \|_\infty = 1 \qquad \Longrightarrow \qquad \beta (C_\phi) = 1 \, .
\end{displaymath} 
The converse implication is true by \eqref{always}. \par
\medskip

The proof follows the same pattern as in \cite{JAT}, with the following additional argument. 

\begin{lemma} \label{enplus} 
Let $\nu$ be a probability  measure, compactly carried by $\phi (\D)$, and let $R_\nu \colon {\mathfrak B} \to L^{2} (\nu)$ be the canonical inclusion. 
Then, we have: 
\begin{displaymath} 
a_{n} (C_\phi) \gtrsim a_{n} (R_\nu) \, .
\end{displaymath} 
\end{lemma}

To prove this lemma, we need another one. For $f \in {\cal H} (\D)$ and $0  < r < 1$, we set as usual:
\begin{displaymath} 
M (r, f) = \sup_{| z | = r} | f (z) | \, . 
\end{displaymath} 
We then have:
\begin{sublemma} \label{huriez} 
Let $g \in {\cal H} (\D)$, not identically zero, and $0 < r < 1$. Then, there exists $C > 0$, depending only on $g$ and $r$, such that: 
\begin{equation} \label{apriori} 
\qquad M (r, f) \leq C \, \| f g \|_{\mathfrak B}, \quad \forall f \in {\cal H} (\D) \, . 
\end{equation}
Therefore, for each compact subset $L \subseteq \D$, there exists a constant $C = C (L, g)$ such that, for any $f \in {\cal H} (\D)$, one has: 
\begin{equation} \label{aposteriori} 
\| f \|_{{\cal C} (L)}\leq C \, \| f g \|_{\mathfrak B} \, .
\end{equation}
\end{sublemma}

\noindent{\bf Proof.} Since the zeros of $g$ are at most countable, we can find $r \leq \rho < 1$ such that $g$ does not vanish on the circle of radius $\rho$. 
Hence there is some $\mu_r > 0$ such that:
\begin{displaymath} 
\qquad | g (a) | \geq \mu_r > 0 \qquad \text{for } | a | = \rho \, .
\end{displaymath} 
Let $\delta = 1 - \rho$, $f \in {\cal H} (\D)$ and $| a | = \rho$. By subharmonicity of $| f g |^2$, we have:
\begin{displaymath} 
\mu_r^{2} | f (a) |^2 \leq | f (a) \, g (a) |^2 \leq \frac{1}{\delta^2} \int_{D (a,\delta)} | f g |^2 \, dA 
\leq \frac{1}{\delta^2} \int_{\D} | f g |^2 \, dA \, , 
\end{displaymath} 
whence $M (\rho, f) \leq C \, \| f g \|_{\mathfrak B}$ with $C = 1 / (\delta \mu_r) $. But $M (r, f) \leq M (\rho, f)$ by the maximum modulus principle, and 
we get \eqref{apriori}. That ends the proof of Sublemma~\ref{huriez}, since if $L \subseteq \overline{D (0, r)}$ and $f \in {\cal H} (\D)$, then 
$\| f \|_{{\cal C} (L)} \leq M (r, f) \leq C \, \| fg \|_{\mathfrak B}$ by the maximum modulus principle again . \qed 
\medskip 

\noindent{\bf Proof of Lemma~\ref{enplus}.} Let ${\rm D} \colon {\mathcal D}_\ast\to {\mathfrak B}$ be the differentiation operator 
(which is a unitary operator by definition of the norms of these spaces). We can write (see \cite{JAT}, proof of Lemma~3.5): $\nu = \phi (\sigma)$, for some 
probability measure $\sigma$ carried by a compact subset $L$ of $\D$. We then have, for any $f \in  {\mathcal D}_\ast$, with help of Sublemma~\ref{huriez} 
applied to the non-zero function $g = \phi '$:
\begin{align*}
\| R_\nu {\rm D} \, f \|_{L^2 (\nu)}^2 
& = \int | f ' |^2 \, d\nu 
= \int_{L} | f ' \circ \phi |^2 \, d\sigma 
\leq \| f ' \circ \phi \|_{{\cal C} (L)}^2 \\
& \leq C^2 \int_{\D} | f ' \circ \phi |^2 \, | \phi ' |^2 \, dA 
= C^2 \, \| C_{\phi} f \|_{ {\mathcal D}_\ast}^2.
\end{align*}
This implies $a_{n} (R_\nu {\rm D}) \leq C\, a_{n} (C_\phi)$, or, equivalently, $a_{n} (R_\nu) \leq C\, a_{n} (C_\phi)$, since ${\rm D}$ is unitary. \qed 
\medskip

Recall now the following lemmas, borrowed from \cite{JAT} (the first one will be used again latter, in Lemma~\ref{steptwo}). 

\begin{lemma} [{\rm \cite{JAT}, Lemma~3.6}] \label{Seville 2} 
For every $r  \in (0, 1)$ there exist $s = s (r) < 1$ and $f = f_r \in H^\infty$ with the following properties:\par \smallskip 
1) $\lim_{r \to 1^-}s(r) = 1$; \par \smallskip 
2) $\Vert f \Vert_\infty \leq 1$; \par \smallskip 
3) $f ((0, r]) = s\, \partial \D$ in a one-to-one way.
\par \smallskip

\noindent Explicitly, one has:
\begin{equation} \label{notation}
\qquad \qquad s = \e^{- \eps \pi / 2} \qquad \text{with} \qquad \eps = \frac{2 \, \pi}{\log \frac{1 + r}{1 - r}} \, \cdot 
\end{equation} 
\end{lemma}
Note that in \cite{JAT}, we defined $\eps$ with the help of a parameter $\rho$, but $\frac{1 + \rho}{1 - \rho} = \sqrt{\frac{1 + r}{1 - r}}$. \par

\begin{lemma} [{\rm see \cite{JAT}, Lemma~3.7 and its proof}] \label{plusplus} 
Let $0 < r < 1$ and $s$ be as in \eqref{notation}. Then, there exists a probability measure $\mu$ carried by $[0, r]$ such that, if  
$R_\mu \colon {\mathfrak B} \to L^{2}(\mu)$ is the canonical inclusion, one has, for every $n \geq 1$: 
\begin{displaymath} 
a_{n} (R_\mu) \gtrsim \frac{s^n}{\sqrt n} \, \cdot
\end{displaymath} 
\end{lemma}
\begin{lemma} [{\rm see \cite{JAT}, Lemma~3.8 and its sequel}] \label{plupluplu} 
\hskip 2pt Let $\phi \colon \D \to \D$ be a Schur function and suppose that $0$ and $r$ belong to $\phi (\D)$. Then, for any probability measure $\mu$ 
carried by $[0, r]$, there exists a probability measure $\nu$ compactly carried by $\phi (\D)$ such that 
\begin{displaymath} 
a_{2n} (R_\mu) \leq 2 \, a_{n} (R_\nu) \, .
\end{displaymath} 
\end{lemma}

\noindent{\bf Proof of Theorem~\ref{seul cas geom} .} The three lemmas put together give the result. Indeed, assume that $\| \phi \|_\infty > r$. By making a 
rotation, we may assume that $r \in \phi (\D)$. Let then $\mu $ be as in Lemma~\ref{plusplus} and $\nu$ be as in Lemma~\ref{plupluplu} (that we may 
use since $0 = \phi (0) \in \phi (\D)$). Using Lemma~\ref{enplus}, we obtain:
\begin{displaymath} 
a_{n} (C_\phi) \gtrsim a_{n} (R_\nu) \gtrsim a_{2n} (R_\mu) \gtrsim \frac{s^{2n}}{\sqrt {2 n}} \raise 1pt \hbox{,}   
\end{displaymath} 
and, changing $s$ into $s^2$, this ends the proof of Theorem~\ref{seul cas geom}. \qed 


\subsection{Extremal behavior} 

In this section, we see that we may have very compact composition operators on ${\mathcal D}_\ast$ whose image touches the boundary of $\D$. \par

\begin{theorem} \label{katkat}  
For every vanishing sequence $(\eps_n)_n$ of positive numbers, there exists a symbol $\phi$ with $\| \phi \|_\infty = 1$ and such that $C_\phi$ is compact on 
${\mathcal D}_\ast$,  but:
\begin{equation} \label{bone} 
a_{n} (C_\phi) \lesssim \e^{- n \eps_n} . 
\end{equation} 
In particular, $C_\phi$ may be in all Schatten classes $S_p ({\mathcal D}_\ast)$, $p > 0$, of the Dirichlet space.
\end{theorem} 

The result will follow from the forthcoming theorem, which is the analogue of Theorem~5.1 in \cite{JAT}. 

\begin{theorem} \label{zinc} 
Let $\phi$ be a Schur function inducing a bounded composition operator on ${\mathcal D}_\ast$.  Set:
\begin{equation} \label{natation} 
m (t) = \frac{1}{t^2} \int_{| w | \geq 1 - t} n_\phi \, d A \qquad \text{and}  \qquad M (t) = \sum_{k = 0}^\infty m (2^{- k} t )  \, .
\end{equation}
Then:
 \begin{equation} \label{generalupper} 
a_{n} (C_\phi) \lesssim \inf_{0 < t < 1} \Big[ n (1 - t)^n + \sqrt{M (t)} \Big] \,.
\end{equation}
\end{theorem}

\noindent {\bf Proof.} We shall need the following simple inequalities. \par\smallskip

1) For $f \in {\mathcal D}_\ast$ and $a \in \D$, one has: 
\begin{equation} \label{reproderiv} 
| f ' (a) | \leq \frac{\| f \|_{\cal D}}{1 - | a |^2} \, \cdot 
\end{equation}
This is clear by using the Cauchy-Schwarz inequality or by using the formula $f ' (a) = \langle f, \frac{\partial{K}}{\partial{\bar{a}}} (a) \rangle$ (where $K$ 
is the reproducing kernel of  $\mathcal{D}_{*}$). \par

2) If $g \in {\mathcal D}_\ast$ and $g = G'$, then:
\begin{equation} \label{shift} 
\| g \|_{\mathcal D} \leq \| z g \|_{\mathcal D} \quad \text{and} \quad \| G \|_{\cal D} \leq \| g \|_{\cal D} \, . 
\end{equation}
This is obvious by inspection of coefficients. \par
\smallskip

Let now $R = C_\phi  S_{n - 1} \colon {\mathcal D}_{*} \to {\mathcal D}_{*}$ be the operator of rank $< n$ defined by:
\begin{displaymath} 
R (f) = \sum_{k = 1}^{n - 1} \widehat{f} (k) \, \phi^k, 
\end{displaymath} 
so that $(C_\phi - R) (f) = C_{\phi} (u)$ with:
\begin{displaymath} 
u (z) = \sum_{k = n}^\infty \widehat{f} (k) \, z^k := z^n \, v (z) 
\end{displaymath} 
and $v \in \mathcal{D}$. \par
Assume once and for all in the proof  that $\| f \|_{\cal D} \leq 1$. Then $\| v \|_{\cal D} \leq  \| u \|_{\cal D} \leq \| f \|_{\cal D} \leq 1$. \par

Fix $0 < h < 1$. We have, writing $v = w '$,  $u ' (z) = n z^{n - 1} w ' (z) + z^n v ' (z)$, and using \eqref{shift}: 
\begin{align*}
\| (C_\phi - R) (f) \|_{\cal D}^2 
& = \| C_\phi (u) \|_{\cal D}^2 
= \int_\D | u ' (z) |^2 \, n_\phi (z) \, d A (z) \\ 
& = \int_{| z | \leq 1 - h} | u ' (z) |^2 \, n_\phi (z) \, d A (z) \\ 
& \hskip 60 pt + \int_{1 - h \leq | z | < 1} | u ' (z) |^2 \, n_\phi (z) \, d A (z) \\ 
& \lesssim n^2 (1 - h)^{2 n} \int_{| z | \leq 1 - h} | w ' (z) |^2 \, n_\phi (z) \, d A (z) \\
& \hskip 50 pt + (1 - h)^{2 n} \, \int_{| z | \leq 1 - h} | v ' (z) |^2 \, n_\phi (z) \, d A (z) \\
& \hskip 50 pt + \int_{1 - h \leq | z | < 1} | u ' (z) |^2 \, n_\phi (z) \, d A (z) \\
& := I_1 + I_2 + I_3.
\end{align*}

Clearly, using that $\| w \|_{\cal D} \leq \| v \|_{\cal D} \leq 1$, we have: 
\begin{displaymath} 
I_1 \leq n^2 (1 - h)^{2 n} \| C_\phi (w) \|_{\cal D}^2 \lesssim n^2 (1 - h)^{2 n} .
\end{displaymath} 
Similarly, we have:
\begin{displaymath} 
I_2 \leq (1 - h)^{2 n} \|  C_\phi (v) \|_{\cal D}^2 \lesssim (1 - h)^{2 n} .
\end{displaymath} 
We estimate $I_3$ by splitting into dyadic annuli: $I_3 = \sum_{k = 0}^\infty J_k$ with:
\begin{displaymath} 
J_k = \int_{1 -2^{- k} h \leq | z | < 1 - 2^{- k - 1} h} | u ' (z) |^2 \, n_\phi (z) \, d A (z) \, .
\end{displaymath} 
We now use the pointwise estimate \eqref{reproderiv} to get: 
\begin{displaymath} 
| u ' (z) |^2 \leq \frac{\| u \|_{\cal D}^2} {(1 - | z |)^2} \leq  \frac{1}{(1 - | z |)^2} \, \cdot 
\end{displaymath} 
In view of \eqref{notation}, this gives an estimate of the form: 
\begin{displaymath}
J_k \leq 4^k h^{- 2} \int_{| z | \geq 1 - 2^{- k} h} \hskip - 5 pt n_\phi (z) \, d A (z) = m (2^{- k} h) \, . 
\end{displaymath}
Summing up, we get $I_3 \leq M (h)$.  It follows that 
\begin{displaymath} 
\| (C_\phi - R) (f) \|^2 \lesssim [n^2 (1 - h)^{2 n} + M (h)]. 
\end{displaymath} 
Taking the supremum on $f$, and then square roots, we then get:
\begin{displaymath} 
a_{n} (C_\phi) \leq  \| C_\phi  - R \| \lesssim [n (1 - h)^n + \sqrt{M (h)}].
\end{displaymath} 
Finally, taking the infimum on $h$, we end the proof of Theorem ~\ref{zinc}. \qed 

\bigskip 

\noindent{\bf Remark.} In \cite{JAT} (Theorem~4.1), we proved in the opposite direction, following \cite{Carroll-Cowen}, that a composition operator on the 
weighted Bergman space ${\mathfrak B}_\alpha$ may be compact, but no little more. It is likely that the same occurs in ${\mathcal D}_\ast$, namely that for 
every vanishing sequence $(\eps_n)_n$ of positive numbers, there exists a symbol $\phi$ such that $C_\phi$ is compact on ${\mathcal D}_\ast$ and for 
which $\liminf_{n \to \infty} \frac{a_{n} (C_\phi)}{\eps_n} > 0$ (in particular, if it happens to be true, we might have $C_\phi$ compact and in no Schatten 
class $S_p ({\mathcal D}_\ast)$, $p < \infty$, of the Dirichlet space). But we do not succeed in proving that.


\section{Approximation numbers of the cusp map}

In \cite{estimates}, it is shown that there is a composition operator $C_\chi \colon H^2 \to H^2$, whose symbol is called the cusp map, defined on the 
Hardy space, such that, for some constants $c_1  > c_2 > 0$, one has:
\begin{equation} 
\qquad \e^{- c_1 n / \log n} \lesssim a_n (C_\chi \colon H^2 \to H^2) \lesssim \e^{- c_2 n / \log n} \, , \qquad n = 2, 3 , \ldots \, .
\end{equation} 

In \cite{Dirichlet}, we proved that every composition operator which is compact on the Dirichlet space is in all Schatten classes $S_p (H^2)$, $p > 0$, on the 
Hardy space. Therefore the approximation numbers of $C_\phi \colon H^2 \to H^2$ must be (much) smaller than those of 
$C_\phi \colon {\mathcal D}_\ast \to {\mathcal D}_\ast$. The next theorem gives, for the cusp map, this order of smallness. 

\begin{theorem} \label{cusp} 
Let $\chi$ be the cusp map. There exist two constants $0 < c' < c$ such that the approximation numbers $a_n (C_\chi)$ of the associated composition operator 
$C_\chi \colon {\cal D}_\ast \to {\cal D}_\ast$ satisfy:
\begin{equation} \label{estimations cusp} 
\qquad \quad \e^{- c \, \sqrt n} \lesssim a_n (C_\chi \colon {\cal D}_\ast \to {\cal D}_\ast) \lesssim \e^{ - c ' \, \sqrt n} \qquad n = 1, 2, \ldots \, .
\end{equation} 
\end{theorem} 

Recall the definition of the cusp map $\chi$, introduced in \cite{LLQR Bergman}, and later used, with a slightly different definition in \cite{estimates}. 
Actually, as in \cite{Dirichlet}, we have to modify it slightly again in order to have $\chi (0) = 0$. We first define:
\begin{displaymath}
\chi_0 (z) = \frac{\displaystyle \Big( \frac{z - i}{i z - 1} \Big)^{1/2} - i} {\displaystyle - i \, \Big( \frac{z - i}{i z - 1} \Big)^{1/2} + 1} \, ; 
\end{displaymath}
we note that $\chi_0 (1) = 0$, $\chi (- 1) = 1$, $\chi_0 (i) = - i$, $\chi_0 (- i) = i$, and $\chi_0 (0) = \sqrt{2} - 1$. Then we set: 
\begin{displaymath}
\chi_1 (z) = \log \chi_0 (z), \quad \chi_2 (z) = - \frac{2}{\pi}\, \chi_1 (z) + 1, \quad \chi_3 (z) = \frac{a}{\chi_2 (z)}  \, \raise 1pt \hbox{,} 
\end{displaymath}
and finally:
\begin{displaymath}
\chi (z) = 1 - \chi_3 (z) \, ,
\end{displaymath}
where:
\begin{equation} \label{definition de a}
a = 1 - \frac{2}{\pi} \log (\sqrt{2} - 1) \in (1, 2) 
\end{equation} 
is chosen in order that $\chi (0) = 0$. The image $\Omega$ of the (univalent) cusp map is formed by the intersection of the  inside of the disk 
$D \big(1 - \frac{a}{2} \raise 1pt \hbox{,} \frac{a}{2} \big)$ and the outside of the two disks  
$D \big(1 + \frac{i a}{2} \raise 1pt \hbox{,} \frac{a}{2} \big)$ and $D \big(1  - \frac{ i a}{2} \raise 1pt \hbox{,} \frac{a}{2} \big)$. 


\subsection{Proof of the upper bound of Theorem~\ref{cusp}} 

We need some lemmas. 

\begin{lemma} 
We have:
\begin{equation} \label{debut} 
\| \chi^n \|_{\cal D} \leq C\, n^{- \delta},
\end{equation} 
for every $n \geq 1$, where $C$ and $\delta$ are positive numerical constants.
\end{lemma} 

\noindent{\bf Proof.} Since $\chi$ is univalent, we have, for every $0 < h < 1$: 
\begin{displaymath} 
\| \chi^n \|_{\cal D}^2 = \int_\D n^2 |w|^{2 n - 2} n_\chi (w) \, d A (w) 
\leq n^2 (1 - h)^{2 n - 2} + n^2 \int_{| w | \geq 1- h} n_\chi (w) d A (w) 
\end{displaymath} 
But $\int_{| w | \geq 1- h} n_\chi (w) \, d A (w)$ is the area of $\chi (\D) \cap \{ |w| \geq  1 - h\}$; since $\chi (\D)$ is delimited at the cuspsidal point $1$ 
by two circular arcs, this area is $\approx h^3$. We get hence:
\begin{displaymath} 
\| \chi^n \|_{\cal D} \lesssim n [ \e^{- n h} + h^{3/2} ] \, .
\end{displaymath} 
The choice $h = 2 (\log n / n)$ gives $\| \chi^n \|_{\cal D} \lesssim n^{- 1 / 2} (\log n)^{3 / 2}$ and hence the lemma, with any $\delta < 1/ 2$. \qed 
\par \medskip 

An immediate corollary, in which $S_N$ denotes the operator of $N$th-partial sum, as defined in \eqref{partial sum}, is the following. 
\begin{corollary} 
We have:  
\begin{displaymath} 
\| C_\chi - C_\chi S_N \| \lesssim N^{- \delta} .
\end{displaymath} 
\end{corollary}

\noindent{\bf Proof.} Using the Hilbert-Schmidt norm, we get: 
\begin{align*}
\hskip 50 pt plus 20 pt minus 20 pt 
\| C_\chi - C_\chi S_N \|^2 
& \leq \| C_\chi - C_\chi S_N \|_{HS}^2 
= \sum_{n > N} \frac{\| \chi^n \|^2}{n} \\ 
& \lesssim \sum_{n > N} n^{- 1 - 2 \delta} \lesssim N^{- 2\delta} \, . \hskip 115 pt plus 20 pt minus 20 pt \square 
\end{align*}

Now, the idea for majorizing $a_n (C_\chi)$ is to write, for every operator $R$ with rank $< n$: 
\begin{displaymath} 
\| C_\chi - R \| \leq \| C_\chi - C_\chi S_N \| + \| C_\chi S_N - R \| \, ,
\end{displaymath} 
which gives, taking the infimum over all such $R$:
\begin{equation} 
a_n (C_\chi) \leq \| C_\chi - C_\chi S_N \| + a_n (C_\chi S_N) \, .
\end{equation} 
Using the corollary, we get:
\begin{equation} 
a_n (C_\chi) \lesssim N^{- \delta} + a_n (C_\chi S_N)
\end{equation} 
and our goal is to give a good upper bound of $a_n (C_\chi S_N)$. \par
\smallskip

\begin{lemma} \label{lundi} 
For some numerical constant $\eps > 0$, we have: 
\begin{equation} \label{optimistic} 
a_n (C_\chi S_N) \lesssim \sqrt{N} \, \e^{- \eps \sqrt{n}}. 
\end{equation}
\end{lemma}
\smallskip 

With this estimation, we get: 
\begin{displaymath} 
a_n (C_\chi) \lesssim [N^{- \delta} + \sqrt{N} \, \e^{- \eps \sqrt{n}}] 
\end{displaymath} 
and, by adjusting $N = \big[ \e^{\eps \sqrt{n}} \big]$, we obtain the upper bound in \eqref{estimations cusp}.  \par
\medskip
 
\noindent{\bf Proof of Lemma~\ref{lundi}.} To prove \eqref{optimistic}, we shall replace $C_\chi S_N$ by a ``dominating'' operator. \par
\smallskip

We begin with observing that, if $f (z) = \sum_{j = 1}^\infty c_j z^j \in {\mathcal D}_{*}$, we have by the change of variable formula, setting 
$d \mu = n_{\chi} \, d A = \ind_{\chi (\D)} \, d A$:
\begin{equation} \label{sentinel} 
\begin{split}
\Vert C_\chi S_N f \Vert_{{\mathcal D}}^2 
& =\int_{\D} \bigg| \sum_{j = 1}^N j \, c_j w^{j - 1}\bigg|^2 n_{\chi} (w) \, d A (w) \\ 
& = \int_{\D} \bigg| \sum_{j = 1}^N j \, c_j w^{j - 1}\bigg|^2 \, d \mu (w) \, .
\end{split}
\end{equation}

Now, denote by $\Delta_N \colon {\mathcal D}_{*} \to H^2$ the map defined by:
\begin{displaymath} 
\Delta_N f (w) = \sum_{j = 1}^N j c_j w^{j - 1} . 
\end{displaymath} 
Observe that: 
\begin{displaymath} 
\| \Delta_N f \|_{H^2}^2 = \sum_{j = 1}^N j^2 | c_j |^2 
\leq N \sum_{j = 1}^N j \, | c_j |^2 
\leq N \| f \|_{\mathcal D}^2,
\end{displaymath} 
so that $\| \Delta_N \| \leq \sqrt N$. \par 

Let also $J$ be the canonical inclusion $J \colon H^2 \to L^2 (\mu)$. The  equality \eqref{sentinel} reads  
$\| C_\chi S_N f \|_{\cal D}^2 = \| J \Delta_N f \|_{L^2 (\mu)}^2$; therefore there is a contraction 
$C_N \colon  {\mathcal D}_{*} \to {\mathcal D}_{*}$ such that
\begin{equation} 
C_\chi S_N = C_N J \Delta_N \, .
\end{equation} 
The ideal property of approximation numbers now implies: 
\begin{displaymath} 
a_n (C_\chi S_N) = a_n (C_N J \Delta_N) \leq \| C_N \| \cdot a_n (J) \cdot \| \Delta_N \| \leq \sqrt{N} \, a_n (J) \, ,
\end{displaymath} 
and we are left with the task of majorizing $a_n (J)$. To that effect, we use the Gelfand numbers $c_n$ (\cite{PIE} or \cite{CAST}). Recall that if 
$T \colon X \to X$ is an operator on some Banach space $X$, then $c_n (T) = \inf \{ \| T_{\mid Z} \| \, ; \  Z \subseteq X, {\rm codim} \, Z < n \}$, and 
if $X = H$ is a Hilbert space, then $c_n (T) = a_n (T)$.\par

Let $B$ be a Blaschke product of length $< n$, let $E = B H^2$ which is a subspace of $H^2$ of codimension $< n$. We have:
\begin{displaymath} 
a_n (J) = c_n (J) \leq \| J_{\mid E} \| \, . 
\end{displaymath} 

The majorization is then made using the Carleson embedding theorem. Let $r$ be the greatest integer $< \sqrt n$, and $B_0$ a Blaschke product with $r$ zeros 
well distributed on the interval $(0, 1)$. More precisely, $B_0$ has its zeros at the points 
\begin{displaymath} 
\qquad z_j = 1 - 2^{- j}, \qquad 1 \leq j \leq r \, .
\end{displaymath} 

Set $\Omega = \chi (\D)$ and observe that: 
\begin{align}
z \in \Omega \text{ and } \Re z \geq 1 - h \quad &\Longrightarrow \qquad \quad  | \Im z | \lesssim h^2 \label{observation} \\
A [S (\xi, h) \cap \Omega] \lesssim h^3 & \qquad \text{for every } \xi \in \T . \label{matage} 
\end{align}

Let now $B = B_{0}^r$. This is a Blaschke product of length $r^2 < n$. Using the Carleson embedding theorem (for the measure $d \mu = n_\chi \, d A$),  
as in \cite{JAT} and \cite{estimates}, and the univalence of $\chi$, we get: 
\begin{equation} \label{carle} 
\| J_{\mid E} \|^2 \lesssim \sup_{0 < h < 1,\ | \xi | = 1}  \frac{1}{h} \int_{S (\xi, h) \cap \, \Omega} | B |^2 \, d A \, .
\end{equation}
To estimate the supremum in the right-hand side of \eqref{carle}, we may assume that $h = 2^{- l}$ and we separate two cases.  \par
\smallskip

$\bullet$ $l \geq r$. Then, using \eqref{matage} and the fact that $| B | \leq 1$, we have: 
\begin{equation} \label{morceau 1} 
\frac{1}{h} \int_{S (\xi, h) \cap \, \Omega} | B |^2 \, d A \lesssim \frac{1}{h}\, h^3 = h^2 = 2^{- 2l} \leq  2^{- 2 r}. 
\end{equation} 

$\bullet$ $l < r$. Then, we have: 
\begin{align*}
\frac{1}{h} \int_{S (\xi, h) \cap \, \Omega} | B |^2 \, d A 
& \leq \frac{1}{h} \int_{\{ | z - 1 | \leq 2^{- r} \} \cap \, \Omega} | B (z) |^2 \, d A (z) \\ 
& +\sum_{j = l + 1}^r \frac{1}{h} \int_{C_j \cap \, \Omega} | B |^2 \, d A \, , 
\end{align*}
where $C_j$ is the annulus
\begin{displaymath} 
C_j = \{ z \in \D \, ; \ 2^{- j} \leq | z - 1 | \leq 2^{- j + 1} \}. 
\end{displaymath} 
The first term is handled as before. Now, since $\Omega$ is contained in some sector $1 - | z | \geq \delta \, | 1 - z |$, 
we have, for $z \in C_j \cap \Omega$:
\begin{displaymath} 
1 - | z | \geq \delta\, 2^{- j} \quad \text{and} \quad  1 - | z_j | = 2^{- j}, 
\end{displaymath} 
whereas 
\begin{displaymath} 
| z - z_j | = | z - 1 + 2^{- j} | \leq | z - 1 | + 2^{- j} \leq 3.2^{- j}.
\end{displaymath} 
This implies that, for some absolute constant $M > 0$:
\begin{displaymath} 
| z - z_j | \leq M \, \min (1 - | z|, 1 - | z_j |) 
\end{displaymath} 
and, by \cite{LLQR Lens}, Lemma~2.3, the $j$-th factor of $B_0$ is, in modulus, less than $\kappa = \frac{M}{\sqrt{M^2 + 1}} < 1$. Therefore 
$| B | = | B_0^r | \leq \kappa^r$ on all sets $C_j \cap \Omega$, so that
\begin{align*}
\sum_{j = l + 1}^r \frac{1}{h} \int_{C_j \cap \, \Omega} | B |^2 \, d A 
& \leq \sum_{j = l + 1}^r \frac{1}{h} \int_{C_j \cap \, \Omega} \kappa^{2 r} \, d A \\ 
& \lesssim 2^l \kappa^{2 r} \, A [ S( \xi, 2^{- l}) ] 
\lesssim  \kappa^{2 r} 2^l 2^{- 2 l} \lesssim \kappa^{2 r} .  
\end{align*}

This finally shows, thanks to \eqref{carle} and \eqref{morceau 1}, that $\| J_{\mid E} \| \lesssim \kappa^r$, or, in setting $\kappa = \e^{- \eps}$, that 
(recall that $r$ is the greatest integer $< \sqrt n$, and hence $r \approx \sqrt n$):
\begin{displaymath} 
\| J_{\mid E} \| \lesssim \e^{-\eps \sqrt{n}}.
\end{displaymath} 

This proves \eqref{optimistic} and ends the proof of  the upper bound in Theorem~\ref{cusp}. \qed 


\subsection{Proof of the lower bound of Theorem~\ref{cusp}} 

Recall that $\mu$ is the measure $d \mu = n_\chi \, d A$ and that $\Omega = \chi (\D)$. \par 

Consider the diagram 
\begin{displaymath} 
\xymatrix{
H^2 \ar[r]^{{\rm P}} & {\mathcal D}_\ast \ar[r]^{C_\chi} & {\mathcal D}_\ast \ar[r]^-{{\rm D}} & L^2 (\mu) \, ,
}
\end{displaymath} 
in which 
\begin{displaymath} 
{\rm P} \,  \Big(\sum_{n = 0}^\infty c_n z^n \Big) = \sum_{n = 0}^\infty c_n \frac{z^{n + 1}}{n + 1} 
\end{displaymath} 
is the ``primitivation'' operator and ${\rm D}$ is the differentiation operator. We have: 
\begin{displaymath} 
{\rm D} \,  C_\chi \, {\rm P}\, f = (f \circ \chi) \, \chi ' \, .  
\end{displaymath} 
We note  that, by definition of the norms, $\| {\rm P} \| \leq 1$. For $0 < h < 1$ fixed, let also: 
\begin{displaymath} 
R \colon H^2 \to L^{\infty} ([0, 1 - h])  
\end{displaymath} 
be the canonical injection. \par

The rest of the proof consists of two steps, the first of which consists of showing that  $a_n (C_\chi)$ is not much smaller than $a_n (R)$. 
\begin{lemma} \label{stepone} 
We have: 
\begin{displaymath} 
a_n (C_\chi) \geq \frac{h^2}{4} \, a_n (R) \, .
\end{displaymath} 
\end{lemma} 

\noindent{\bf Proof.} We first notice that, if $f \in H^2$, and $0 \leq x \leq 1 - h$, we have: 
\begin{equation} \label{etun} 
\| R (f) \|_{L^{\infty}([0, 1 - h])} \leq \frac{4}{h^2} \| f \|_{L^2 (\mu)}. 
\end{equation}
To that effect, we observe that (recall that $a \in (1, 2)$ in defined in \eqref{definition de a}):
\begin{equation} \label{vu} 
0 < h \leq a - 1 \quad \text{and} \quad 0 \leq x \leq 1 - h \quad \Longrightarrow \quad D (x, h^2 / 4 a) \subseteq \Omega \, .
\end{equation}
Indeed, if $z = u + i v \in D (x, h^2 /4 a)$ and $0 \leq x \leq 1 - h$, we have $1 - u \geq h - (h^2 / 4 a) \geq h/ 2$, as well as $| v | < h^2 / 4 a$,  and:
\begin{displaymath} 
\Big| z - \Big(1 + \frac{i a}{2} \Big) \Big|^2 = (1 - u)^2  + \Big(v - \frac{a}{2} \Big)^2 
\geq \frac{h^2}{4} + v^2 - a \, |v| + \frac{a^2}{4} > \frac{a^2}{4} \, \cdot
\end{displaymath} 
Similarly $\big| z - \big(1 - \frac{i a}{2} \big) \big| > \frac{a}{2}$. Moreover, since $1 - \frac{a}{2} \leq \frac{a}{2} - h$, we have  
$\big| z - \big( 1 - \frac{a}{2} \big) \big| \leq | z - x|+ \big|x - \big( 1 - \frac{a}{2} \big) \big| \leq \frac{h^2}{4a}+\frac{a}{2} - h < \frac{a}{2}$. 
Hence $z \in \Omega$. \par

Therefore, by subharmonicity of the function $| f |^2$:
\begin{align*} 
| f (x) |^2 
& \leq \frac{16 \, a^2}{h^4} \int_{D (x, h^2 / 4)} | f |^2 \, d A \\ 
& \leq  \frac{16 \, a^2}{h^4} \int_{\Omega} | f |^2 \, d A 
= \frac{16 \, a^2}{h^4} \int_{\D} | f |^2 n_\chi \, d A \\ 
& = \frac{16 \, a^2}{h^4} \int_{\D} | f  |^2 \, d \mu \, ,
\end{align*} 
which proves \eqref{etun}. 
\par\medskip 

Let now $f \in H^2$ and $g = {\rm P} \, f \in {\mathcal D}_{*}$, so that $f = {\rm D} \,  g$.  As follows from \eqref{etun} and from the change of variable 
formula, we have: 
\begin{align*}
\| R f \|_{\infty}^2 
& \leq \frac{64}{h^4} \int_{\D} | f (w) |^2 n_\chi (w) \, d A (w)  
= \frac{64}{h^4} \int_{\D} | {\rm D} \,  g (w) |^2 n_\chi (w) \, d A (w) \\ 
& = \frac{64}{h^4} \int_{\D} | g' \big(\chi (z) \big) |^2 | \chi ' (z) |^2 \, d A (z)  
=\frac{64}{h^4} \| {\rm D}\,  C_\chi g \|_{L^2 (\D)}^2 \\ 
& = \frac{64}{h^4} \| C_\chi {\rm P}  f \|_{\mathcal{D}}^2.  
\end{align*}
Therefore, there exists $C \colon \mathcal{D}_{*} \to L^{\infty} ([0, 1 - h])$ such that:
\begin{displaymath} 
R = C \, C_\chi {\rm P} \qquad \text{and} \qquad \| C \| \leq \frac{4}{h^2} \, \cdot
\end{displaymath} 
All this implies, by the ideal property of approximation numbers: 
\begin{displaymath} 
a_n (R) \leq \| C \| \,  a_n (C_\chi) \, \| {\rm P} \| \leq \frac{4}{h^2} \, a_n (C_\chi) \, , 
\end{displaymath} 
which ends the proof of Lemma~\ref{stepone}. \qed 

\par\medskip 

The second step  consists of a minoration of $a_n (R)$, which uses the comparison with Bernstein numbers and a good choice of an $n$-dimensional space $E$. 

\begin{lemma} \label{steptwo} 
Let $ 0 < r <1$ and $s$ as in \eqref{notation}. We have: 
\begin{displaymath} 
a_n (R) \geq  \frac{s^n}{\sqrt n} \, \cdot 
\end{displaymath} 
\end{lemma} 

Recall (see \cite{PIE} for example) that, if $X$ and $Y$ are two Banach spaces, and $T \colon X \to Y$ is a compact operator, the $n$-th Bernstein number 
of $T$ is: 
\begin{displaymath} 
b_n (T) = \sup_{{\rm dim}\,  E = n} \inf_{f \in S_E} \| T f \| \, , 
\end{displaymath} 
where $S_E$ denotes the unit sphere of $E$, and we have:
\begin{equation}  \label{beber} 
a_n (T) \geq b_n (T) \, .
\end{equation} 

\noindent{\bf Proof of Lemma~\ref{steptwo}.} Let $f = f_r$ be as in Lemma~\ref{Seville 2}, and write $r = 1 - h$. Consider the $n$-dimensional space 
\begin{displaymath} 
E = [1, f, \ldots, f^{n - 1} ]  \, ,  
\end{displaymath} 
and let $g = \sum_{j = 0}^{n - 1} \alpha_j f^j \in E$ with $\| g \|_\infty = 1$. By Lemma~\ref{Seville 2} and the Cauchy-Schwarz inequality, we have:
\begin{displaymath} 
1 \leq \sum_{j = 0}^{n - 1} | \alpha_j | \,  \| f^j \|_\infty 
\leq \sum_{j = 0}^{n - 1} | \alpha_j |  
\leq \sqrt{n} \Big( \sum_{j = 0}^{n - 1} | \alpha_j |^2 \Big)^{1 / 2} .
\end{displaymath} 
On the other hand, Lemma~\ref{Seville 2} again gives us: 
\begin{align*}
\| R (g) \|_\infty 
& \geq \Big\| \sum_{j = 0}^{n - 1} \alpha_j \, s^j \, \e^{i j \theta} \Big\|_{L^\infty (\T)} 
\geq \Big\| \sum_{j = 0}^{n - 1} \alpha_j \, s^j \, \e^{i j \theta} \Big\|_{L^2 (\T)} \\ 
& = \Big( \sum_{j = 0}^{n - 1} | \alpha_j |^2 \, s^{2 j} \Big)^{1 / 2} 
\geq s^n \Big( \sum_{j = 0}^{n - 1} | \alpha_j |^2 \Big)^{1 / 2} \\ 
& \geq \frac{s^n}{\sqrt n} \, \cdot
\end{align*}
Therefore, $b_n (R) \geq s^n / \sqrt{n}$. Using \eqref{beber}, we get  $a_n (R) \geq s^n / \sqrt{n}$ as well.  \qed 

\par \bigskip 

Let us now indicate how Lemma~\ref{stepone} and Lemma~\ref{steptwo} allow to finish the proof. Write $h = \e^{- A}$ where $A > 0$. Then, with the 
notation \eqref{notation}, we have:
\begin{displaymath} 
\eps = \frac{2 \, \pi} {\log \frac{1 + r}{1 - r}} \lesssim \frac{1}{\log \frac{1}{1 - r}} = \frac{1}{\log (1 / h)} = \frac{1}{A} 
\end{displaymath} 
and 
\begin{displaymath} 
s \gtrsim \e^{- c / A} \, ,
\end{displaymath} 
for some constant $c > 0$. Therefore  Lemma~\ref{stepone} and Lemma~\ref{steptwo} give: 
\begin{displaymath} 
a_n (C_\chi) \gtrsim h^2 a_n (R) \gtrsim h^2 \frac{s^n}{\sqrt n} \gtrsim \e^{- c ' \, (A + n/ A)}. 
\end{displaymath} 

The optimal choice $A = \sqrt n$ gives the lower bound in Theorem~\ref{cusp}. \qed 

\bigskip

\noindent{\bf Remark:} One sees that the approximation numbers of $C_\chi$ behave quite differently on the Hardy space $H^2$ (like 
$\e^{- c \, n / \log n}$, see \cite{estimates}) and on the Dirichlet space (like $\e^{- c \, \sqrt n}$). This seems to be due to the following. 
On the Hardy space, the important fact is the parametrization $t \mapsto \chi (\e^{i t})$ where logarithms are involved. On the Dirichlet space, 
we only need to know the geometry of $\chi (\D)$, a domain limited by three circles, where logarithms are no longer involved. 

\goodbreak


\section{ Capacity of the set of contact points}

Here is now the improvement of  a theorem in (\cite{ELFKELSHAYOU}) in terms of approximation numbers (see also \cite{Dirichlet}, Theorem~4.1). 
This improvement is definitely optimal in view of our previous Theorem~\ref{seul cas geom}, stating that,  for every bounded composition operator $C_\phi$ 
on ${\mathcal D}_\ast$ of symbol $\phi$, one has: 
\begin{displaymath} 
\| \phi \|_\infty = 1 \qquad \Longrightarrow \qquad \beta (C_\phi) := \liminf_{n \to \infty} [ (a_n (C_\phi) ]^{1 / n} = 1.
\end{displaymath} 

Recall the following notation, where $\phi$ belongs to the disk algebra $A (\D)$, {\sl i.e.} the space of continuous functions $f \colon \overline{\D} \to \C$ which are analytic in $\D$: \par

\begin{displaymath} 
E_\phi = \{ \e^{i t} \in \T \, ; \ |\phi (\e^{i t} )| = 1 \} \, .
\end{displaymath} 
\begin{theorem} \label{improvement} 
Let $K$ be a compact set of the circle $\T$ with  logarithmic capacity ${\rm Cap}\, K = 0$, and $(\eps_n)_n$ a sequence of positive numbers with limit $0$. 
Then, there exists a Schur function $\phi$ generating a composition operator $C_\phi$ bounded on ${\mathcal D}$ and with the following properties:  \par
\smallskip 
1)  $\phi \in A (\D) \cap {\mathcal D}$, the ``Dirichlet algebra''; \par
\smallskip
2) $E_\phi = K$ and $E_\phi = \{ \e^{it} \in \T \, ; \, \phi (\e^{it}) = 1 \}$; \par
\smallskip
3) $a_n (C_\phi) \lesssim \e^{- n \,\eps_n}$ for all $ n \geq 1$.
\end{theorem}

Before proving this theorem, we need two results \cite{Dirichlet}. The first one is the existence of a peculiar peaking function. Recall that a function 
$q \in A (\D)$ is said to peak on a compact subset $K \subseteq \partial{\D}$, and is called a peaking function, if:
$q (z) = 1$ for $z\in K$ and $| q(z) | < 1$ for $z \in \overline{\D} \setminus K$. \par

\begin{theorem} [\cite{Dirichlet} Theorem~4.2] \label{q peak}
For every compact set $K \subseteq \partial{\D}$ of logarithmic capacity ${\rm Cap} \, K = 0$, there exists a Schur function $q \in A (\D) \cap {\cal D}_\ast$ 
which peaks on $K$ and such that the composition operator $C_q \colon {\cal D}_\ast \to {\cal D}_\ast$ is bounded (and even Hilbert-Schmidt).
\end{theorem}

The other one is a lemma borrowed from the proof of Theorem~3.3 in \cite{Dirichlet}.  

\begin{lemma} \label{recall} 
Let $\delta \colon (0,1) \to (0, \infty)$ be a positive function with $\lim_{h \to 0} \delta (h) = 0$. Then, there exists a univalent Schur function 
$\gamma \in A (\D)$ such that $\gamma (1) = 1$  and that:
\begin{equation} \label{small} 
\int_{| w | \geq 1 - h} n_{\gamma} (w) \, d A (w) = A [\gamma (\D) \cap \{ w \, ; \ 1 - h \leq | w | < 1 \}] \leq \delta (h) \, . 
\end{equation}
\end{lemma}

\noindent{\bf Proof of Theorem~\ref{improvement}.} It suffices to use Lemma~\ref{recall} to construct a generalized cusp map $\gamma$ in order to have 
$a_n (C_\gamma) \lesssim \e^{- n \eps_n}$. Indeed, we then compose this generalized cusp map $\gamma$ with a symbol $q$ peaking on $K$, as given by 
Theorem~\ref{q peak}; namely consider $\phi = \gamma \circ q$. Then, we know that $E_\phi = \{ \e^{i t} \, ; \ \phi (\e^{i t}) = 1 \} = K$. Moreover,  
$C_\phi = C_q \circ C_\gamma$, so that, using the fact that $\| C_q \| < \infty$: 
\begin{displaymath} 
a_n (C_\phi) \leq \| C_q \| \,  a_n (C_\gamma) \lesssim \e^{- n \eps_n} \, .
\end{displaymath} 

It remains to find such a generalized cusp map $\gamma$. Set: 
\begin{displaymath} 
\delta_n = \eps_n + \frac{\log n}{n} \, \cdot
\end{displaymath} 

Let $\Phi$ be a positive, continuous, concave, and increasing function $ \Phi \colon [0, 1] \to [0, 1]$ such that $\Phi (0) = 0$ and 
$\Phi (1 / n) \geq \delta_n$. Let $\Psi = \Phi^{- 1}$. By Lemma~\ref{recall}, we can adjust $\gamma$ so as  to have, using the notation \eqref{natation}: 
\begin{displaymath} 
m (h) \leq h \, \rho^2 (h) \, ,
\end{displaymath} 
where  
\begin{displaymath} 
\rho (h) = \exp \Big( - \frac{h}{\Psi (h)} \Big) \, \cdot
\end{displaymath} 
Note that $\rho$ is increasing. We then see that: 
\begin{displaymath} 
M (h) \leq \sum_{k = 0}^\infty 2^{- k} h \, \rho^2 (2^{- k} h) \lesssim h\, \rho^2 (h) \leq \rho^2 (h)
\end{displaymath} 
and, plugging that in \eqref{generalupper}, we get: 
\begin{displaymath} 
a_n (C_\gamma) \lesssim \inf_{0 < h < 1} [n (1 - h)^n + \rho (h) ] 
\leq \inf_{0 < h < 1} \bigg[n \exp(- n h) + \exp \Big(- \frac{h}{\Psi(h)} \Big) \bigg] \, .
\end{displaymath} 

In particular, if we choose $h = \Phi (1 / n)$, we obtain: 
\begin{displaymath} 
a_n (C_\gamma) \lesssim n \, \e^{- n\Phi (1 / n)} \leq n \, \e^{- n \delta_n} = \e^{- n \eps_n} .
\end{displaymath} 

In view of the initial observation, this ends the proof of Theorem~\ref{improvement}. \qed



\bigskip

\vbox{\small 
\noindent{\it 
{\rm Pascal Lef\`evre}, Univ Lille-Nord-de-France UArtois, \\ 
Laboratoire de Math\'ematiques de Lens EA~2462, \\
F\'ed\'eration CNRS Nord-Pas-de-Calais FR~2956, \\
F-62\kern 1mm 300 LENS, FRANCE \\ 
pascal.lefevre@euler.univ-artois.fr 
\smallskip

\noindent
{\rm Daniel Li}, Univ Lille-Nord-de-France UArtois, \\
Laboratoire de Math\'ematiques de Lens EA~2462, \\
F\'ed\'eration CNRS Nord-Pas-de-Calais FR~2956, \\
Facult\'e des Sciences Jean Perrin,\\
Rue Jean Souvraz, S.P.\kern 1mm 18, \\
F-62\kern 1mm 300 LENS, FRANCE \\ 
daniel.li@euler.univ-artois.fr
\smallskip

\noindent
{\rm Herv\'e Queff\'elec}, Univ Lille-Nord-de-France USTL, \\
Laboratoire Paul Painlev\'e U.M.R. CNRS 8524, \\
F-59\kern 1mm 655 VILLENEUVE D'ASCQ Cedex, FRANCE \\ 
queff@math.univ-lille1.fr
\smallskip

\noindent
{\rm Luis Rodr{\'\i}guez-Piazza}, Universidad de Sevilla, \\
Facultad de Matem\'aticas, Departamento de An\'alisis Matem\'atico \& IMUS,\\ 
Apartado de Correos 1160,\\
41\kern 1mm 080 SEVILLA, SPAIN \\ 
piazza@us.es\par}
}

\end{document}